\documentclass{commat}

\title{Transcendental Continued Fractions}

\author{Sarra Ahallal, Ali Kacha}

\affiliation{
    \address{Department of Mathematics, Faculty of Science,
     Ibn Tofail University, 14 000 Kenitra, Morocco}
\email{sarraahallal92@gmail.com, ali.kacha@uit.ac.ma}
}

\abstract{%
    In the present paper, we give sufficient conditions on the elements of the continued fractions $A$ and $B$ that will assure us that the continued fraction $A^B$ is a transcendental number. With the same condition, we establish a transcendental measure of  $A^B.$
   }



\VOLUME{30}
\NUMBER{1}
\firstpage{251}
\DOI{https://doi.org/10.46298/cm.10519}

\begin{paper}

\section{Introduction and preliminaries}
The theory of transcendental numbers has a long history.
We know since J. Liouville in $1844$ that the very rapidly converging sequences
of rational numbers provide examples of transcendental numbers. So, in his famous paper
 \cite{p7} Liouville showed that a real number admitting very
good rational approximation could not be algebraic, and then he explicitly constructed
the first examples of transcendental numbers.

From this sense, the transcendence of the continued
fractions having partial quotients that increase rapidly have been
studied by several authors such as P. Bundschuh \cite{ 1}, A. Durand
\cite{p2}, W. Lianxiang \cite{p6}, G.Nettler \cite{p8}, T. Okano \cite{p10}.  We also note that the transcendence of
some power series with rational coefficients is given by some authors, see \cite{p5},
\cite{p9}.

 Let $A$ and $B$ be
two continued fractions which are defined by
$$
A= a_{0}+\frac{1|}{|a_{1}}+ \frac{1|}{|a_{2}}+\cdot \cdot \cdot
$$
and
 $$
 B=b_{0}+\frac{1|}{|b_{1}}+ \frac{1|}{|b_{2}}+\cdot \cdot
\cdot,
$$
 where $(a_{i}>0$, $b_{i}>0 $ are integers for any $i\geq 1).$

In \cite{p3} we have proved the transcendence of the six numbers $%
A $, $B$, $A\pm B$, and $AB^{\pm 1\text{ }}$ if
$a_{n}>b_{n}>a_{n-1}^{\alpha}$ where $\alpha$ is a real constant
$>7$ by using Roth's approximation theorem.
The algebraic independence of $A$ and $B$ was also proved under a similar method in \cite{p4}.

 We recall that in \cite{p1}, P. Bundschuh
has noted that obviously no transcendence proof of $A^B$ could be established by using Roth's approximation
theorem. Then he gave a theorem which yields the transcendence of $A^B$ by using another method.

 In the present paper, the first aim is to improve Bundschuh's result concerning the transcendence of $A^B$
with a slight hypothesis.

 The second main result of this article is to establish a
transcendental measure of a continued fraction $A^B.$

 In order to prove the transcendence of continued fraction $A^B$, we will use the
result of Bundshuh which expresses a necessary condition if $A^B$ is an algebraic number.

  \begin{theorem}[\cite{p11}]
 Let $\xi $ be a real number, $\delta$ a real number $>2,$ if there
exists an infinity rational numbers $\frac{p}{q}$ with $gcd(p,q)=1$ such that
\begin{equation*}
\left\vert \xi -\frac{p}{q}\right\vert <\frac{1}{q^{\delta}},
\end{equation*}
then $\xi$ is a transcendental number.
\end{theorem}

\section{Main results}
\subsection{Transcendence}

With the same notations as above for  $A$ and $B,$ we put
\begin{eqnarray*}
A_{n} =a_{0}+\frac{1|}{|a_{1}}+ \frac{1|}{|a_{2}}+\cdot \cdot+ \frac{1|}{|a_{n}}=\frac{%
^{a}p_{n}}{^{a}q_{n}}, \\
B_{n} = b_{0}+\frac{1|}{|b_{1}}+ \frac{1|}{|b_{2}}+\cdot \cdot+ \frac{1|}{|b_{n}} =\frac{%
^{b}p_{n}}{^{b}q_{n}}.
\end{eqnarray*}

 The principal result of the transcendence of $A^B$ is given in the following  theorem.

\begin{theorem} \label{th 1}
    Let $(a_n)$, $(b_n)$ and $A$, $B$ be as before and let $1 < \alpha \le \alpha'$ be real numbers. If $a_{n+1} \geq b_{n+1} \geq a_{n}^{\alpha },$
 $ a_n^{\alpha'} \geq a_{n+1}$ for all $n \geq 1$ and
 $ \ln(a_n) \ln(b_{n })=O( \ln(b_{n+1}\;^{b}q_{n}^{2} ))$ then $ A^B$ is a
transcendental number.
\end{theorem}

  \begin{remark}
  \begin{enumerate}
      \item We note that under the hypothesis
$a_{n+1} \geq b_{n+1} \geq a_{n}^{\alpha }$ with $\alpha >1,$  the real numbers $A$ and $B$ are transcendental, see \cite{p3}.

\item In \cite{p8}, Nettler presented certain explicit formulae for the algebraic numbers
$A_n^{B_n}$ which converges to $A^B$ as $n \rightarrow + \infty,$
depending on the $m^{th}$ convergents $(m=1,\cdot \cdot\cdot,n)$ of $A$ and $B.$
\end{enumerate}
 \end{remark}

  The proof of the Theorem \ref{th 1} is based on some lemmas that
express the link between partial quotients of the mentioned reals
continued fractions and the denominators of their convergents.

 The main tool of this Theorem is the following necessary condition for algebricity
of $A^B$ due to Bundschuh.

 \begin{lemma}\label{lem 1}
 Let $A,$ $ B$ as before, but such that $A^B$ is an algebraic number and that the sequence $(^{b}q_{n})_n$
 satisfies
 \begin{equation}\label{1}
 \ln (^{b}q_{n+1})=O(^{b}q_{n} \ln (^{b}q_{n}))
 \end{equation}
 as  $n \rightarrow + \infty.$ Then there exists an effectively computable number $\gamma_0 >0,$ depending only on $A$
 and $B$ such that the inequality

 \begin{equation}\label{2}
 \max(|A-\beta|, |B-\delta|) \leq \exp(-\gamma_0 \ln( H_1) \ln(H_2) )
 \end{equation}
 has at most finitely many solutions $(\beta, \delta,H_1,H_2) \in \mathbb{Q}^2\times \mathbb{N}^2$
 with $h(\beta) \leq H_1,$ $ h(\delta) \leq H_2,$ and $H_1,H_2 \geq 4,$  $h(\theta) $ denotes the height of an algebraic number $\theta.$
 \end{lemma}

\begin{lemma}\label{lem 2}
\begin{enumerate}
\item  Let $\alpha >1$ be a real number. If  $a_{n} \geq a_{n-1}^{\alpha }$ for all $n \geq 2,$ then we have
 $$
 ^{a}q_{n}<  2 a_{n}^{\frac{\alpha }{\alpha -1}}.
 $$
 \item  Let $\alpha >1$ be a real number. If $a_{n} \geq b_{n} \geq a_{n-1}^{\alpha}$  for all $ n \geq 2,$ then there exists a real constant $C_1>0$ such that
\begin{equation}\label{3}
^{a}q_{n}  >  {}\text{ }^{b}q_{n} > C_1 {}\text{
}^{a}q_{n-1}^{\alpha }.
\end{equation}
\end{enumerate}\end{lemma}

\begin{remark} The same result as (1) of Lemma \ref{lem 2} is also obtained for $^{b}q_{n}$ and $b_n.$ \end{remark}

 \begin{proof}[Proof of Lemma \ref{lem 2}] \begin{enumerate}
     \item  See \cite{p3}.

  \item We can easily show that
\begin{eqnarray}\label{4}
^{b}q_{n} &=&b_{n}\text{ }^{b}q_{n-1}+\text{ }^{b}q_{n-2} \nonumber \\
& > &b_{n}\text{ }^{b}q_{n-1}\geq a_{n-1}^{\alpha }\text{
}^{b}q_{n-1\text{
\quad }}  \nonumber \\
& > &\prod\limits_{i=1}^{n-1}a_{i}^{\alpha }.
\end{eqnarray}
 On the other hand, one has
\begin{equation*}
^{a}q_{n-1} < (a_{n-1}+1)^{a}q_{n-2} < \prod\limits_{i=1}^{n-1}\left( 1+\frac{1}{a_{i}}\right)
\prod\limits_{i=1}^{n}a_{i},
\end{equation*}
 which becomes
\begin{equation}\label{5}
\prod\limits_{i=1}^{n-1}a_{i}^{\alpha }> ^{a}q_{n-1}^ {\alpha} \prod\limits_{i=1}^{n-1}
\left( 1+\frac{1}{a_{i}}\right)^{-\alpha}.
\end{equation}
 Combining $(4)$ and $(5)$ gives
\begin{equation*}
^{b}q_{n} > \text{ }^{a}q_{n-1}^{\alpha }
\prod\limits_{i=1}^{n-1}\left( 1+\frac{1}{a_{i}}\right)^{-\alpha}> C_2 \text{
}^{a}q_{n-1}^{\alpha }
\end{equation*}
 where $C_2=C_2(A,\alpha)>0.$
We notice that the left side of the inequality $(3)$ can be proved in a similar way.
\qedhere
\end{enumerate}
\end{proof}

\begin{proof}[Proof of Theorem \ref{th 1}] First, let us justify that the hypothesis $(1)$ of Lemma \ref{lem 1} is verified.
According to Lemma \ref{lem 2}, for any $\varepsilon>0$, we have
$$
^{b}q_{n+1}<b_{n+1}^{\frac{\alpha}{\alpha-1}+\epsilon}<b_{n}^{\frac{\alpha.\alpha'}{\alpha-1}+\varepsilon \alpha'}
$$
 for any sufficiently large $n.$ Which yields
\begin{align*}
	\ln(^{b}q_{n+1})&<(\dfrac{\alpha.\alpha'}{\alpha-1}+\varepsilon \alpha')\ln(b_{n})\\
	&<(\dfrac{\alpha.\alpha'}{\alpha-1}+\varepsilon \alpha')\ln(^{b}q_{n}) \\
	&<(\dfrac{\alpha.\alpha'}{\alpha-1}+\varepsilon \alpha')\ln(^{b}q_{n})^{b}q_{n}.
	\end{align*}
 So, the relationship $(1)$ is verified, it suffices to make $\varepsilon$ tend to $0.$

 We will show that under the hypothesis of Theorem \ref{th 1} the relationship $(2)$ is satisfied by the infinite numbers:
$$(\dfrac{^{a}p_{n}}{^{a}q_{n}} \;\;, \dfrac{^{b}p_{n}}{^{b}q_{n}}\;\;,H_{1}=\;^{a}q_{n}=h(A_{n})\;\; ,H_{2}=\;^{b}q_{n}=h(B_{n}))\;\in\; \mathbb{Q}^{2}\times\mathbb{N}^{2}.
$$

 By part $2$ of Lemma \ref{lem 2} we have $ ^{b}q_{n} < ^{a}q_{n},$ so
 \begin{equation}\label{6}
max(\mid A-A_{n} \mid,\mid B-B_{n} \mid )<\dfrac{1}{^{b}q_{n+1}\;^{b}q_{n}}<\dfrac{1}{b_{n+1}\;^{b}q_{n}^{2}}
\end{equation}
 for all $n \geq 1.$

 On the other hand, the result $1$ of Lemma \ref{lem 1} implies that

\begin{equation}\label{7}
\ln(^{a}q_{n}).\ln(^{b}q_{n})<(\dfrac{\alpha}{\alpha - 1})^{2}\ln(a_{n}).\ln(b_{n}).
\end{equation}

 Namely, the hypothesis of Theorem \ref{th 2}
 $$\ln(a_{n}).\ln(b_{n})=O(\ln(b_{n+1}\;^{b}q_{n}^{2})),$$
 shows that there exists a real constant $C_2 >0,$ such that
\begin{equation}\label{8}
	\ln(a_{n}).\ln(b_{n}) < C_{2}\ln(b_{n+1}\;^{b}q_{n}^{2})
\end{equation}	
  for all sufficiently large $n.$ Combining $(7)$ and $(8),$ we get
\begin{equation}\label{9}
\ln(^{a}q_{n}).\ln(^{b}q_{n}) < C_{3}\ln(b_{n+1}\;^{b}q_{n}^{2}),
\; C_{3}=(\dfrac{\alpha}{\alpha-1})^{2}C_{2}
\end{equation}
 for any $n$ sufficiently large.

 Due to inequality \eqref{9}, the relationship \eqref{6} becomes
\begin{align*}
	max(\mid A-A_{n} \mid,\mid B-B_{n} \mid)&<\exp(-\dfrac{1}{C_{3}}\ln(^{a}q_{n})\ln(^{b}q_{n})) \\
	&<\exp(-\dfrac{1}{C_{3}}\ln(h(A_{n}))\ln(h(B_{n})))
\end{align*}
 for infinitely many $n.$ This shows that the conclusion of Lemma \ref{lem 1} fails, hence $A^{B}$
is a transcendental number. \end{proof}

\begin{example}\label{ex1} Let
 $\begin{cases}
	a_0=b_0=0, \; a_{1}=b_{1}=3, \; \delta \;\text{\;a\; real\; number} >0, \\
	a_{n+1}=b_{n+1}=a_{n}^{1+\delta}\;, n\geq1. \\
    \alpha=\alpha'=1+\delta.
\end{cases}
$
\end{example}

  We verify that the hypothesis $\ln(a_{n})\ln(b_{n})=O(\ln(b_{n+1}\;^{b}q_{n}^{2}))$
of Theorem \ref{th 1} is satisfied. We have
$$\ln(a_{n})\ln(b_{n})=(1+\delta)^{2n-2}\times(\ln a_{1})^{2},$$
 and for all $n \geq 2,$ we get
\begin{align*}
b_{n+1}\;^{b}q_{n}^{2}&\geq b_{n+1} b_{n}^{2}\\
&\geq a_{1}^{(1+\delta)^n} \times \big(a_{1}^{(1+\delta)^{n-1}}\big)^2\\
&\geq \big(a_{1}\big)^{(1+\delta)^n+(1+\delta)^{2n-2}}\\
&\geq \big(a_{1}\big)^{(1+\delta)^{2n-2}}.
\end{align*}
 It follows that
$$\dfrac{\ln(a_{n}).\ln(b_{n})}{\ln(b_{n+1}\;^{b}q_{n}^{2})}\leq \ln(a_{1})=\ln 3.
$$
By applying Theorem \ref{th 1}, we deduce that $A^A$ is a transcendental number.

\subsection{The transcendental measure of a continued
fraction $A^B$}

 In this subsection, we give the second main result of
this article. We keep the same notations as in the first subsection.

 \begin{theorem}\label{th 2} Let $P \in {\mathbb{Z}}[X]\ \{0\} $ be a polynomial of degree $d \geq 2,$ and height
$H\geq a_{2}^{d+1/2}.$   Let $k$ be a real number $ \geq 1$ and $\alpha =2d+1$ such that
\begin{equation*}
a_{n}^{\alpha  } \leq b_{n+1} \leq a_{n+1} \leq a_{n}^{k\alpha}
\text{ \ for all \ } n\geq 1.
\end{equation*}
 Then, we have
\begin{equation*}
\left\vert P( A^B) \right\vert > \frac{1}{2} ( Hd(d+1))^{-2kd(d+1) }.
\end{equation*} \end{theorem}

 \begin{remark} We note that the continued fraction $A^B$ defined in Theorem \ref{th 2} is a transcendental number. \end{remark}

 For the proof of Theorem \ref{th 2} further Lemma is needed.

   \begin{lemma}\label{lem 3}  Let $k$ be a real number $>1.$
  The hypothesis $ \ a_{n} \leq a_{n-1}^{k\alpha  }$ for
all $n\geq 2$  implies that%
\begin{equation*}
^{a}q_{n} \ < \ 2 \ ^{a}q_{n-1}^{\frac{k\alpha^2}{\alpha -1}}.
\end{equation*}
\end{lemma}
 \begin{proof} By the same method as in Lemma \ref{lem 2} we
prove Lemma \ref{lem 3}.  \end{proof}

 \begin{proof}[Proof of Theorem \ref{th 2}]
 We have assumed that the two numbers
$A$ and $ B$ are larger than $1,$ so $A_{n}<1$ and $B_{n}<1$ for all $n\geq 1.$ Put
$$
P(X)=\sum_{k=0}^{d}e_{k}X^k, \; e_{k}\in \mathbb{Z}.
$$
 In order to obtain a transcendental measure of $A^B,$ we must  minus $\mid P(A^B)\mid$ by a positive function of $H$ and $d.$

 From equality
$$P(A_{n}^{B_{n}})=P(A^{B})+P(A_{n}^{B_{n}})-P(A^B),
$$
 we get
\begin{align*} \mid P(A_{n}^{B_{n}})\mid \leq\mid P(A^{B})\mid +\mid P(A_{n}^{B_{n}})-P(A^B)\mid.
\end{align*}
 Therefore,
\begin{align}\label{10}
\mid P(A^B) \mid\geq \mid P(A_{n}^{B_{n}})\mid - \mid P(A_{n}^{B_{n}})-P(A^B)\mid.
\end{align}

 Firstly, we have
 \begin{align*}
	P(A_{n}^{B_{n}})&=\sum_{k=0}^{d}e_{k}\dfrac{(^{a}p_{n})^{kB_{n}}}{(^{a}q_{n})^{kB_{n}}}\\
&=\dfrac{1}{^{a}q_{n}^{dB_{n}}}\sum_{k=0}^{d}e_{k} (^{a}p_{n})^{kB_{n}}(^{a}q_{n})^{(d-k)B_{n}}.
\end{align*}

  Notice that
  $$\sum_{k=0}^{d}e_{k}(^{a}p_{n})^{kB_{n}} (^{a}q_{n})^{(d-k)B_{n}} \neq 0,
  $$
   because, if we assume that
  $$
  \sum_{k=0}^{d}e_{k}(^{a}p_{n})^{kB_{n}} (^{a}q_{n})^{(d-k)B_{n}} =  (^{a}q_{n})^{(d-k)B_{n}}.P(A_{n}^{B_{n}} =0,
  $$
   then we would have
  $$P(A^B)=lim_{n\rightarrow +\infty}P(A_{n}^{B_{n}})=0$$
   which contradicts the transcendence of $A^B.$ It follows that
\begin{equation}\label{11}
\mid P(A_{n}^{B_{n}}) \mid \geq \dfrac{1}{^{a}q_{n}^{dB_{n}}} > \dfrac{1}{^{a}q_{n}^{d}}.
\end{equation}

 On the other hand, according to Roll's theorem applied to $P$ in the interval $[A^B, A_{n}^{B_{n}}]$ or $[A_{n}^{B_{n}}, A^B],$
 there exists a real number $E\in ]A^{B};A_{n}^{B_{n}}[$ or $ ]A_{n}^{B_{n}},A^B[$ such that:
	\begin{equation}\label{12}
	 \mid P(A^B)-P(A_{n}^{B_{n}})\mid = \mid A^{B}-A_{n}^{B_{n}} \mid  \mid P'(E) \mid, \text{  \; with \; }
0<E<1.
\end{equation}
 From this, we obtain
\begin{equation}\label{13}
\mid P'(E) \mid \leq \sum_{k=1}^{d}k \mid e_{k}	\mid \;\leq H\dfrac{d(d+1)}{2},
\end{equation}
and the relationship $\eqref{12}$ becomes
 \begin{equation}\label{14}
 \mid P(A^B)-P(A_{n}^{B_{n}}) \mid \leq \dfrac{1}{2}Hd(d+1)\mid A^{B}-A_{n}^{B_{n}}\mid.
 \end{equation}
 We notice that the function $ f(x,y)= x^y$ is continuously differentiable on every compact subset $K$ of
$[0,1] \times [0,1].$ Then by the mean value theorem to $f,$ it follows the existence of a real constant $C_1=C(A,B)>0$ depending only on $A$ and $B$
such that
\begin{equation*}
\mid A^B-A_{n}^{B_{n}} \mid  \leq C_{1} \max(\mid A-A_{n}\mid,\mid B-B_{n}\mid)
	                          < C_{1}\mid B-B_{n} \mid  < \dfrac{C_{1}}{b_{n+1}},
\end{equation*}
 for sufficiently large $n.$ Therefore, one has
\begin{equation*}
	\mid P(A^B)-P(A_{n}^{B_{n}}) \mid <C_{1}H\dfrac{d(d+1)}{2}\dfrac{1}{b_{n+1}} <\dfrac{C_{1}Hd(d+1)}{2a_{n}^{\alpha}},
\end{equation*}
 for sufficiently large $n.$

 Using relationships $\eqref{11}$ and $\eqref{14},$  $\eqref{10}$ becomes
\begin{eqnarray}\label{15}
\mid P(A^B) \mid &> \dfrac{1}{^{a}q_{n}^{d}}-\dfrac{C_{1}Hd(d+1)}{2a_{n}^{\alpha}}
      &>\dfrac{1}{{a_n}^{d}} -\dfrac{C_{1}Hd(d+1)}{2a_{n}^{\alpha}}\cdot
\end{eqnarray}
 Now, we look for an integer $n$ from which the quantities $\dfrac{1}{a_{n}^d} $
and $ \dfrac{C_{1}Hd(d+1)}{ 2a_{n}^{\alpha}}$
are of the same order. So, in order to satisfy
 $
 \mid P(A^B) \mid > 1/(2a_{n}^d),
 $
  it is sufficient to have
  $$\dfrac{1}{a_{n}^{d}} - \dfrac{C_{1}Hd(d+1)}{2a_n^{\alpha}} > \frac{1}{2 a_{n}^{d}} \cdot
  $$
 Which is equivalent to the following inequality
 $$\dfrac{1}{2a_{n}^{d}}> \dfrac{C_{1}Hd(d+1)}{2a_n^{\alpha}}$$
  that is
 $$
 a_n^{\alpha} >  C_{1}Hd(d+1)a_{n}^{d}.
 $$
 For this, it is sufficient that $n$ verifies
\begin{equation}\label{16}
 a_{n}^{\frac{\alpha}{2}} > C_{1} \;
	\text{ \; and \; } \;
 a_{n}^{\frac{\alpha}{2}} > Hd(d+1)a_{n}^{d}.
\end{equation}
 The first part of Relationship $\eqref{16}$ is easy to achieve. Therefore, the second inequality  is equivalent to
$$
 a_{n}^{\frac{\alpha}{2}-d}>Hd(d+1).
 $$
 Let $n_{1}$ be the smallest integer $\geq 2$ such that
\begin{equation}\label{17}
a_{n_1-1}^{\frac{\alpha}{2}-d} \leq Hd(d+1) < a_{n_1}^{\frac{\alpha}{2}-d} \cdot
\end{equation}

 The integer $n_{1}$ exists because we have assumed  that $Hd(d+1)>H\geq\; a_2^{\alpha/2}$
and the sequence $(a_{n}^{\alpha})_n$ is increasing and tends to $+\infty$. Then, we have
$$
\mid P(A^{B})\mid  > \dfrac{1}{2a_{n_1}^{d}} \cdot
$$
 Since, $ a_{n_{1}} \leq a_{n_{1}-1}^{k\alpha},$ then we obtain
$$\mid P(A^{B})\mid  > \dfrac{1}{{2a_{n_{1}-1}}^{kd \alpha}} \cdot
$$
On the other hand, the left-hand side of $\eqref{17}$ implies that
\begin{equation*}
	\mid P(A^B)\mid
	>\dfrac{1}{2(HD(d+1))^{\frac{kd\alpha}{\alpha/2-d}}} \cdot
\end{equation*}
 Finally, since $\alpha=2d+1,$ we conclude that
\begin{equation*}
	\mid P(A^B)\mid
	>\dfrac{1}{2}(HD(d+1))^{-2kd(2d+1)} \cdot
\end{equation*}
 Which achieves the proof of Theorem \ref{th 2}. \end{proof}

 \begin{example}
    Let
\[
\begin{cases}
    a_{0}=b_{0}=0, a_{1}=b_{1}=2, \\
	a_{n+1}=b_{n+1}=a_n^5,\;n \geq 2, \\
	\alpha=\alpha'=5, k=1.
\end{cases}
\]
 \end{example}
 In the same way, as in Example \ref{ex1}, we prove that $A^A$ is a transcendental number.
 Let $P\in \mathbb{Z}[X]\setminus\{0\}$ be a quadratic polynomial of height $H\geq\;a_{2}^{5/2}=2^{25/2}.$
By applying Theorem \ref{th 2}, a transcendental measure of $A^A$ is given by
\begin{equation*}
\mid P(A^A)\mid > \frac{1}{2} (6H)^{-10}.
\end{equation*}

\EditInfo{October 19, 2021}{February 21, 2022}{Attila Bérczes}

\end{paper}